\documentclass[groupedaddress,notitlepage,10pt,preprintnumbers]{revtex4-1}
\usepackage{amsmath, amssymb}
\usepackage{cleveref}
\usepackage{csl}
\usepackage{ifthen}
\usepackage{textcomp}
\usepackage{tabularx}
\usepackage{listings}
\usepackage{graphicx}
\usepackage[dvipsnames]{xcolor}

\graphicspath{{figures/}}

\def\!#1{\mathcal{#1}}

\newcommand{\R}{\mathbb{R}}
\newcommand{\C}{\mathbb{C}}
\newcommand{\dd}{\mathrm{d}}
\newcommand{\ddt}[2][]{\frac{\mathrm{d}^{#1} #2}{\mathrm{d} t^{#1}}}

\newcommand{\MATLAB}{\texttt{MATLAB}}

\newcommand{\problemdescription}[3]{
\subsection{#1}
\noindent\textbf{Path:} \texttt{csl.otp.#2}\\
}

\gdef\internalotpname{}
\newenvironment{presettable}[1]{\noindent\tabularx{\textwidth}{|l|X|c|c|c|}
\hline%
\textbf{Name} & \textbf{Path} & \multicolumn{1}{|l|}{\textbf{\textnumero Vars}} & \multicolumn{1}{|l|}{\textbf{Stiff}} & \multicolumn{1}{|l|}{\textbf{Chaotic}}\gdef\internalotpname {#1}\\\hline%
}{\endtabularx}

\newcommand{\preset}[6]{
\textbf{#1} & \texttt{csl.otp.\internalotpname{}.presets.#2} &
#3 &
\ifthenelse{\equal{#4}{Yes}}{\checkmark}{\text{\sffamily X}} &
\ifthenelse{\equal{#5}{Yes}}{\checkmark}{\text{\sffamily X}}\\\hline
\multicolumn{5}{|p{0.99\textwidth}|}{#6}\\\hline
}

\newcommand{\featurestable}[5]{
\begin{tabular}{ccccc}
     Jacobian &  Jvp & Javp & Stiff & Chaotic\\
     \ifthenelse{\equal{#1}{y}}{\checkmark}{\text{\sffamily X}} &
     \ifthenelse{\equal{#2}{y}}{\checkmark}{\text{\sffamily X}} &
     \ifthenelse{\equal{#3}{y}}{\checkmark}{\text{\sffamily X}} &
     \ifthenelse{\equal{#4}{y}}{\checkmark}{\text{\sffamily X}} &
     \ifthenelse{\equal{#5}{y}}{\checkmark}{\text{\sffamily X}}
\end{tabular}
}

\newcommand{\getpresettable}[1]{
\ifthenelse{\equal{#1}{bouncingball}}{\begin{presettable}{bouncingball}
\preset{Canonical}{Canonical}{4}{No}{No}{This simply calls the Flat Terrain example.}
\preset{Parabola}{Parabola}{4}{No}{No}{The ball bounces on a parabola, starting slightly off-center as to create
   an interesting trajectory.}
\preset{Random Terrain}{RandomTerrain}{}{No}{No}{}
\preset{Flat Terrain}{Simple}{4}{No}{No}{The ball has no horizontal velocity, and bounces on a perfectly flat surface.}
\end{presettable}
}{}
\ifthenelse{\equal{#1}{brusselator}}{\begin{presettable}{brusselator}
\preset{Decay}{Decay}{2}{No}{No}{Rapid descent to a fixed point}
\preset{Periodic}{Periodic}{2}{No}{No}{A periodic cycle.}
\preset{Spiral}{Spiral}{2}{No}{No}{Rapid decay into a fixed orbit}
\end{presettable}
}{}
\ifthenelse{\equal{#1}{grayscott}}{\begin{presettable}{grayscott}
\preset{Random}{Random}{5,000}{Yes}{No}{Random initial conditions, useful for stiff integrator testing}
\end{presettable}
}{}
\ifthenelse{\equal{#1}{hires}}{\begin{presettable}{hires}
\preset{Canonical}{Canonical}{8}{Yes}{No}{The canonical initial conditions as used in the literature}
\end{presettable}
}{}
\ifthenelse{\equal{#1}{linear}}{\begin{presettable}{linear}
\preset{Dahlquist}{Dahlquist}{1}{No}{No}{The Dahlquist scalar linear test problem.}
\end{presettable}
}{}
\ifthenelse{\equal{#1}{lorenz63}}{\begin{presettable}{lorenz63}
\preset{Canonical}{Canonical}{3}{No}{Yes}{This is the original problem presented in the literature.
   The initial condition is purposefully outside of the trapping region, but converges to it quite quickly.}
\preset{Limit Cycle}{LimitCycle}{3}{No}{No}{A non-chaotic set of parameters for the problem to showcase potential periodic behavior.}
\preset{Surprise}{Surprise}{3}{No}{No}{Strogatz's `surprise'}
\end{presettable}
}{}
\ifthenelse{\equal{#1}{lorenz96}}{\begin{presettable}{lorenz96}
\preset{Canonical}{Canonical}{40}{No}{Yes}{This is the original problem presented in the literature.
   The initial condition is a slight perturbation of a critical point.}
\preset{Time-dependent Forcing}{PopovSandu}{40}{No}{Yes}{Used in (Popov and Sandu, 2019)}
\end{presettable}
}{}
\ifthenelse{\equal{#1}{qgso}}{\begin{presettable}{qgso}
\preset{Gaspari-Cohn}{GC}{16,129}{No}{Yes}{Uses the Gaspari-Cohn compactly supported function to build an
   interesting initial condition}
\end{presettable}
}{}
}

\begin{document}

\title{ODE Test Problems: a MATLAB suite of initial value problems}

\author{Steven Roberts}
\email[]{steven94@vt.edu}
\affiliation{Computational Science Laboratory, Department of Computer Science, Virginia Tech}
\author{Andrey A. Popov}
\email[]{apopov@vt.edu}
\affiliation{Computational Science Laboratory, Department of Computer Science, Virginia Tech}
\date{\today}
\author{Adrian Sandu}
\email[]{sandu@cs.vt.edu}
\affiliation{Computational Science Laboratory, Department of Computer Science, Virginia Tech}

\begin{abstract}
    ODE Test Problems (OTP) is an object-oriented MATLAB package offering a broad range of initial value problems which can be used to test numerical methods such as time integration methods and data assimilation (DA) methods.  It includes problems that are linear and nonlinear, homogeneous and nonhomogeneous, autonomous and nonautonomous, scalar and high-dimensional, stiff and nonstiff, and chaotic and nonchaotic.  Many are real-world problems from fields such as chemistry, astrophysics, meteorology, and electrical engineering.  OTP also supports partitioned ODEs for testing IMEX methods, multirate methods, and other multimethods.  Functions for plotting solutions and creating movies are available for all problems, and exact solutions are provided when available. OTP is desgined for ease of use---meaning that working with and modifying problems is simple and intuitive.
\end{abstract}

\cslauthor{Steven Roberts, Andrey A. Popov, and Adrian Sandu}
\cslrevision{1}
\csltitlepage

\maketitle

\tableofcontents

\section{Introduction}

The initial value problem \cite{hairer1993solving, hairer1996solving} is very often an important sub-problem in many different applications. 
The aim of this project is to provide not only a useful collection of initial value problems, but an easy-to-use interface for working with them. 
Unlike other projects of this type, our aim is to have setup and use be as simple as possible, with almost no boilerplate code, but with the extendability that comes with more sophisticated, in-depth tools. The project was born out of necessity: disjoint collections of problems with their own setup scripts and uses existed internally. Their standardization in a unified framework written in \MATLAB\ is simply the natural evolution, and a sign of things to come. All problems are implemented purely in \MATLAB\ without any external dependencies or \texttt{mex} files.

Using OTP is as simple as
\begin{lstlisting}
model = csl.otp.lorenz96.presets.Canonical;
[t, y] = ode45(model.RHS.F, model.TimeSpan, model.Y0);
model.movie(t, y);
\end{lstlisting}
A non-trivial example showcasing the versatility of the framework is
\begin{lstlisting}
model = csl.otp.qgso.presets.GC('small');
model.TimeSpan = [0, 5000];
[~, y1] = ode45(model.RHS.F, model.TimeSpan, model.Y0);
model.Parameters.linearsolver = 'multigrid';
model.Parameters.linearsolvertol = 1e-8;
[~, y2] = ode45(model.RHS.F, model.TimeSpan, model.Y0);
disp(norm(y1(end, :) - y2(end, :))/norm(y1(end, :)));
   0.821952174839369
\end{lstlisting}
This showcases that not only is it trivial to use complex models, but that modifying fundamental properties of the model is as easy as typing a few lines of code.

\section{Background}

We consider the classic initial value problem
\begin{align}
    \begin{split}
    	\ddt{\*y(t)} &= \*f(t, \*y(t)), \qquad t \in [t_0, t_f] \\
    	\*y(t_0) &= \*y_0,
    \end{split}
\end{align}
where $\*f : \Omega \subset \R \times \C^n \to \C^n$. With all possible trajectories $(\mathfrak{t},\mathfrak{y})$ defining the space $\Omega$.

The three things that uniquely define an initial value problem are the function $\*f$, commonly referred to as the `right hand side' or RHS for short, the time span of the problem (from an initial $t_0$ to a final $t_f$), and the initial value, $\*y_0$, defined to be the known value of the solution at the initial time. 

Additional things that are helpful but not necessary in the solution of an initial value problem numerically might be the knowledge of higher order derivatives of $\*f$, such as the Jacobian, spatial relations between the variables in terms of a distance function, splittings of the RHS into physical processes that occur at different temporal scales, etc.

\section{ODE Test Problems API}

\subsection{Presets}
By far, the most useful part of OTP are presets. A preset provides problem parameters, the time span, and initial conditions in order to simplify the initialization of a model. The calling convention of a preset could not be simpler:
\begin{lstlisting}
model = csl.otp.{problemname}.presets.{Presetname};
\end{lstlisting}
Because of \MATLAB's default support of tab-completion, it is trivial to cycle through all possible problems and all possible presets. 
All presets are callable without any arguments, with sensible defaults being applied, and all problems have a default preset named \texttt{Canonical}, which represents the problem specification that most commonly appears in literature, or a default designed by us. 
Arguments to all presets are passed with Name-Value pairs, the documentation for which appears whenever the \texttt{help} command is issued on the preset:
\begin{lstlisting}
help csl.otp.{problemname}.presets.{Presetname}
\end{lstlisting}

\subsection{Usage}
As before, the three things required for defining an initial value problem are the right hand side function, time span, and initial condition.
These are very easy to access in OTP:
\begin{lstlisting}
% define the model
model = csl.otp.lorenz96.presets.Canonical;
% accessing the right hand side---a function of time and state
model.RHS.F;
% accessing the timespan---a column vector of initial time and final time
model.TimeSpan;
% accessing the initial condition---a column vector
model.Y0;
% getting the value of the right hand side function at the initial condition
model.RHS.F(model.TimeSpan(1), model.Y0);
\end{lstlisting}
Accessing additional properties of the model (if it has them defined) is just as easy:
\begin{lstlisting}
% defining a random column vector of the size of the model state
rv = randn(model.NumVars, 1);
% evaluating the adjoint model at the initial condition on the random vector
model.RHS.JacobianAdjointVectorProduct(model.TimeSpan(1), model.Y0, rv);
\end{lstlisting}
Finally, a practical example. In order to efficiently solve the problem with an implicit time integration method like \texttt{ode15s} built into \MATLAB, one needs to pass the Jacobian to the time integrator:
\begin{lstlisting}
options = odeset('Jacobian', model.RHS.Jacobian);
[t, y] = ode15s(model.RHS.F, model.TimeSpan, model.Y0, options);
model.movie(t, y);
\end{lstlisting}

\subsection{Parameter Validation}
An important design principle of OTP is to protect the user. 
If the user does something outside of the design parameters of the problem, like inputting in impossible initial conditions or impossible parameter values, OTP is smart enough to error out, and give useful errors.
\begin{lstlisting}
model = csl.otp.lorenz63.presets.Canonical;
model.Parameters.rho = -1;
    Error using csl.odeutils.StructParser/checkField (line 26)
    The field rho does not satisfy nonnegative
    
model.Parameters.rho = [1, 1];
    Error using csl.odeutils.StructParser/checkField (line 26)
    The field rho does not satisfy scalar
\end{lstlisting}

\subsection{Python Interface}
\texttt{Python} has an interface with \MATLAB. This interface does amount to writing matlab code in python, but it does make it possible to interface.
\begin{lstlisting}[language=python]
>>> import matlab.engine
>>> me = matlab.engine.start_matlab()
>>> model = me.csl.otp.lorenz63.presets.Canonical()
>>> me.workspace["model"] = model
>>> rval = me.eval("model.RHS.F(model.TimeSpan(1), model.Y0)")
>>> print(rval[0][0])
10.0
\end{lstlisting}

\subsection{Complex Example}
We will compute the Lyapunov exponents of the Lorenz '63 system. We will be using a large range of the APIs found in OTP. 
\begin{lstlisting}
model = csl.otp.lorenz63.presets.Canonical;
% evolve the system initially
[~, y] = ode45(model.RHS.F, model.TimeSpan, model.Y0);
model.Y0 = y(end, :).';
% make the next timespan start at the end of the previous and the 
% span itself be 500 time units
model.TimeSpan = model.TimeSpan(end) + [0; 500];
% set the tolerances to the lowest allowable by MATLAB and solve
options = odeset('AbsTol', 100*eps, 'RelTol', 100*eps);
[t, y] = ode45(model.RHS.F, model.TimeSpan, model.Y0, options);
dts = diff(t);
Q = eye(model.NumVars);
les = zeros(model.NumVars, 1);
J = model.RHS.Jacobian(t(1), y(1, :).');
for i = 1:(numel(dts) - 1)
    % compute the full linearlization of 
    dt = dts(i); 
    % Do a trapezpidal rule step
    W1 = J*Q;
    J = model.RHS.Jacobian(t(i + 1), y(i + 1, :).');
    W2 = J*Q;
    % we want our QR to be in canonical form.
    [Q, R] = qr(Q + (dt/2)*(W1 + W2)); D = diag(sign(diag(R))); Q = Q*D; R = D*R; 
    rd = diag(R); les = les + log(rd);
end
les = les/(t(end) - t(1));
k = sum(arrayfun(@(ki) sum(les(1:ki)), 1:numel(les)) > 0);
fracd = k + sum(les(1:k))/abs(les(k+1));
disp(les);
    0.9102
    0.0027
  -14.5891
  
disp(fracd);
    2.0626
\end{lstlisting}
This is really close to the best known empirical results~\cite{viswanath1998lyapunov}!

\section{List of Test problems}
This section is a comprehensive list of and user guide to all the test problems that are included in the OTP test suite. 

Each listing gives the path in \MATLAB, the corresponding number of variables that the problem contains, and a short description of all the presets with some of their properties.

\problemdescription{Linear}{linear}{$n$ (1)}

One of the simplest and most fundamental ODEs is the linear, first order ODE:
\begin{equation}
    \label{eq:linear_ode}
    \ddt{y} = A(t) \*y.
\end{equation}
It is well known that the exact solution is $\*y(t) = \exp{\left( \int_{t_0}^t A(\tau) \, \dd{\tau} \right)} \*y_0$.  When \cref{eq:linear_ode} is a scalar, constant coefficient ODE, the problem is often called to Dahlquist test problem and is commonly used to access the stability of ODE solvers~\cite{dahlquist1963special}.  

\getpresettable{linear}
\problemdescription{Bouncing Ball}{bouncingball}{4}

The bouncing ball test problem~\cite{tufillaro1992experimental} is meant to test event function support in time integration methods. It describes the mechanics of an inelastic collision in a perfect conservative system. The ball is modeled as a point-mass in a friction-less 2D environment.

The effect of a point mass falling to the Earth can be described by the second order ODE
\begin{align}
    \begin{split}
        \ddt[2]{x} &= 0, \\
        \ddt[2]{y} &= -g,
    \end{split}
\end{align}
where $x$ is the horizontal spacial dimension, $y$ is the vertical spacial dimension, and $g$ is the acceleration due to gravity.  We can rewrite this into a 4 variable first-order ODE as follows:
\begin{align*}
    \ddt{x_1} &= x_2,\\
    \ddt{y_1} &= y_2,\\
    \ddt{x_2} &= 0,\\
    \ddt{y_2} &= -g.
\end{align*}

It is immediately evident that this is a simple linear problem. The uniqueness of this problem comes from the event function
\begin{align*}
    e(x_1,y_1,x_2,y_2) = y_1 - h(x_1),
\end{align*}
where $h$ is the function defining the height of the ground at the current horizontal spatial location.
When the event function is triggered, the time integration method terminates, and the velocities are transformed by
\begin{align*}
    \begin{bmatrix}x_1 \\ y_1\end{bmatrix}_{\text{new}} = {\begin{bmatrix} \cos(\theta) & \sin(\theta)\\ -\sin(\theta) & \cos(\theta)\end{bmatrix}}^{-1}\begin{bmatrix} 1 & 0\\ 0 & -1\end{bmatrix}\begin{bmatrix} \cos(\theta) & \sin(\theta)\\ -\sin(\theta) & \cos(\theta)\end{bmatrix}\begin{bmatrix}x_1 \\ y_1\end{bmatrix},
\end{align*}
where $\theta = \arctan(h'(x_2))$ represents the angle of the ground at the current horizontal position. This boils down to rotating the velocities into a reference frame where the velocity is purely vertical, inverting it, and rotating back into the original reference frame.

This, of course, can be rewritten in the efficient form,
\begin{align*}
\begin{bmatrix}x_1 \\ y_1\end{bmatrix}_{\text{new}} = 
\begin{bmatrix}
    \frac{\left(1 - {h'(x_2)}^2\right) x_1 + 2 h'(x_2) y_1}{1 + {h'(x_2)}^2}\\
    \frac{2 h'(x_2) x_1 - \left(1 - {h'(x_2)}^2\right) y_1 + }{1 + {h'(x_2)}^2}
 \end{bmatrix}.
\end{align*}

An immediate condition as a result of this is that the ground function is differentiable everywhere, and that there are no vertical walls anywhere in the space.

An example of the problem in action is 
\begin{lstlisting}
model = csl.odetestproblems.bouncingball.presets.RandomTerrain;
options = odeset('AbsTol', 1e-6, 'RelTol', 1e-6, 'MaxStep', 2^-4);
[t, y] = csl.odeutils.solveeventproblem(@ode45, model, options);
model.phaseSpacePlot(t, y);
\end{lstlisting}

\begin{figure}
\begin{center}
  \includegraphics[width=0.65\textwidth]{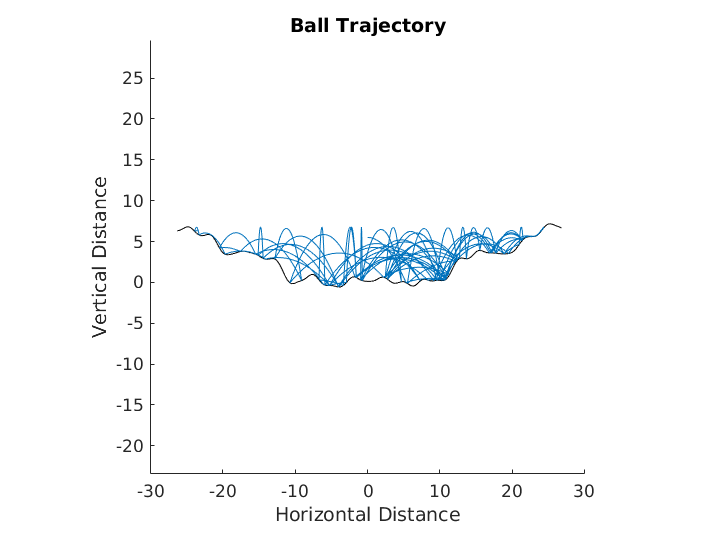}%
\end{center}
\caption{Modeling a perfect bouncing ball in a vacuum, on randomly generated sinusoidal terrain.}
\end{figure}

\getpresettable{bouncingball}
\problemdescription{Brusselator}{brusselator}{2}

The Brusselator is a model for an autocatalytic reaction characterized by the reactions
\begin{align*}
    A & \rightarrow X, \\
    2X + Y & \rightarrow 3X, \\
    B + X & \rightarrow Y + C, \\
    X & \rightarrow D.
\end{align*}
This chemical reaction corresponds to the dimensionless system of first order differential equations
\begin{align}
    \begin{split}
        \ddt{x} &= 1 - (b + 1) x + a x^2 y, \\
        \ddt{y} &= b x - a x^2 y,
    \end{split}
\end{align}
where $a,b$ are positive, real constants and $x,y$ are the concentrations of the two reactants~\cite{ault2003dynamics}.

An interesting splitting of this problem is
\begin{equation*}
    f =
    \underbrace{\begin{bmatrix} 1 - (b + 1) x\\ b x \end{bmatrix}}_{\text{linear}}
    + \underbrace{\begin{bmatrix} a x^2 y\\ -a x^2 y \end{bmatrix}}_{\text{nonlinear}}.
\end{equation*}
This is an additive splitting of the right hand side into linear and nonlinear parts, which can be used in partitioned or exponential time integrators.

\getpresettable{brusselator}
\problemdescription{Double Pendulum}{doublependulum}{4}

The double pendulum problem models the trajectory of two pendulums connected end to end by massless rods.  The first pendulum has a rod of length $ l_{1} $ and a bob of mass $ m_{1} $.  The second pendulum is connected to the first pendulum's bob and has rod length $ l_{2} $ and bob mass $ m_{2} $.  This dynamical system is very sensitive to initial conditions and is governed by
\begin{align}
    \begin{split}
        \ddt[2]{\theta_{1}} &= \frac{m_{2}l_{1}\theta_{1}^2\sin{\Delta}\cos{\Delta}+m_{2}g\sin{\theta_{2}}\cos{\Delta}+m_{2}l_{2}\theta_{2}^2\sin{\Delta}-mg\sin{\theta_{1}}}{ml_{1}-m_{2}l_{1}\cos^2{\Delta}}, \\
        \ddt[2]{\theta_{2}} &= \frac{-m_{2}l_{2}\theta_{2}^2\sin{\Delta}\cos{\Delta}+m(g\sin{\theta_{1}}\cos{\Delta}-l_{1}\theta_{1}^2\sin{\Delta}-g\sin{\theta_{2}})}{ml_{2}-m_{2}l_{1}\cos^2{\Delta}},
    \end{split}
\end{align}
where $g$ is the acceleration due to gravity, $m = m_{1} + m_{2}$, and $\Delta = \theta_{2} - \theta_{1}$.  The angles $\theta_1, \theta_2$ are measured counterclockwise from the negative y-axis.

\getpresettable{doublependulum}
\problemdescription{HIRES}{hires}{8}

This HIRES (High Irradiance Responses) problem is a mildly stiff system of eight first order differential equations.  It corresponds to a chemical reaction that models how light affects morphogenesis in a plant.  The equations are given by
\begin{align}
    \begin{split}
        \ddt{y_1} &= -1.71 y_1 + 0.43 y_2 + 8.32 y_3 + 0.0007, \\
        \ddt{y_2} &= 1.71 y_1 - 8.75 y_2, \\
        \ddt{y_3} &= -10.03 y_3 + 0.43 y_4 + 0.035 y_5, \\
        \ddt{y_4} &= 8.32 y_2 + 1.71 y_3 - 1.12 y_4, \\
        \ddt{y_5} &= -1.745 y_5 + 0.43 y_6 + 0.43 y_7, \\
        \ddt{y_6} &= -280 y_6 y_8 + 0.69 y_4 + 1.71 y_5 - 0.43 y_6 + 0.69 y_7, \\
        \ddt{y_7} &= 280 y_6 y_8 - 1.81 y_7, \\
        \ddt{y_8} &= -280 y_6 y_8 + 1.81 y_7,
    \end{split}
\end{align}
with $\*y_{0}=[1,0,0,0,0,0,0,0.0057]^\intercal$ as the initial condition~\cite{schafer1975new}.

\getpresettable{hires}
\problemdescription{Lorenz '63 (3-variable)}{lorenz63}{3}

The three variable Lorenz problem~\cite{lorenz1963deterministic,strogatz2018nonlinear} (named here as Lorenz '63 to contrast it with the $N$-variable Lorenz '96) is the first concrete ODE described to encompass chaotic dynamics.
\begin{align}
    \begin{split}
        \ddt{x} &= \sigma (y-x),\\
        \ddt{y} &= x(\rho-z)-y,\\
        \ddt{z} &= x y-\beta z,
    \end{split}
\end{align} 
with the canonical initial conditions outside of the trapping region: $[x_0, y_0, z_0] = [0,1,0]^\intercal$

Lorenz '63 is a chaotic problem and, as such, is sensitive to initial conditions and parameters. In this example, we will propagate forward the same initial conditions through a problem with one of the parameters being slightly perturbed.
\begin{lstlisting}
model = csl.otp.lorenz63.presets.Canonical;
[t1, y1] = ode45(model.RHS.F, model.TimeSpan, model.Y0);
pert = sqrt(eps(model.Parameters.sigma));
model.Parameters.sigma = model.Parameters.sigma + pert;
[~, y2] = ode45(model.F, t1, model.Y0);
model.plot(t1, y1 - y2);
\end{lstlisting}
\begin{figure}
\begin{center}
  \includegraphics[width=0.90\textwidth]{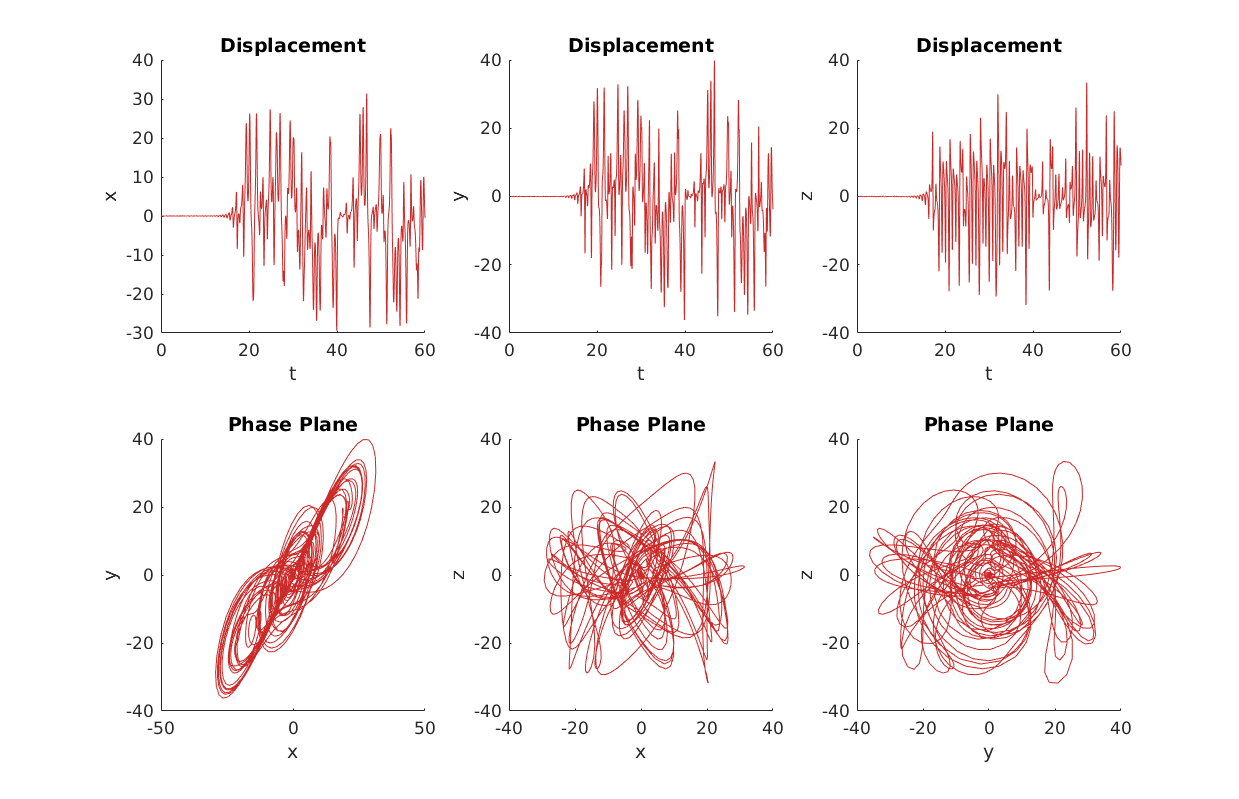}%
\end{center}
\caption{Perturbations in model runs of the canonical Lorenz '63 test problem resulting from a tiny perturbation in the parameter $\sigma$.}
\end{figure}

\getpresettable{lorenz63}
\problemdescription{Lorenz '96}{lorenz96}{$n$ (40)}

The Lorenz '96 system~\cite{lorenz1996predictability} can be defined as follows
\begin{equation*}
f_i = \ddt{x_i} = (x_{i+1} - x_{i-2})x_{i-1} - x_i + F(t),
\end{equation*}
where $i = 1,2,\dots,N$ with $x_{N+1} = x_1, x_{0} = x_N, x_{-1} = x_{N-1}$ and where $F(t)$ is a forcing function usually defined as $F(t) = 8$.
The canonical amount of variables of the system is $N=40$.

In this discretization, one time unit corresponds to 5 days. The canonical time interval for solution is 0.05 units which corresponds to 6 hours.

Lorenz '96 is one of the simplest $N$-dimensional chaotic test problems. In the canonical ($F(t) =8$, $N=40$) case, the system has 13 positive Lyapunov exponents with a fractal dimension of about $27.1$~\citep{popov2018bayesian}.

Note that $x_i = 8$ is a critical point, and is thus likely to be inside a trapping region. A slight perturbation from that is therefore a good choice of initial condition. Canonically, the 40 variable case,
\begin{equation*}
    x_i = \begin{cases} 8.008 & i = 20\\ 8 & \mathfrak{sonst}.\end{cases}
\end{equation*}

\getpresettable{lorenz96}
\problemdescription{Quasi-Geostrophic 1.5 Layer Model (QGSO)}{qgso}{$n^2$ (16,129)}

The quasi-geostrophic model of Sakov and Oke (QGSO)~\cite{sakov2008deterministic} is defined by the PDE
\begin{align}
\begin{split}
\label{eq:QGSO}
\frac{\partial q}{\partial t} = -\psi_x-\varepsilon J(\psi,q)-A\Delta^3\psi+2\pi\sin(2\pi y),\\
\Delta\psi - F\psi = q,\quad J(\psi,q)\equiv\psi_x q_y - \psi_y q_x,
\end{split}
\end{align}
by non-dimensionalizing potential vorticity.
\Cref{eq:QGSO} is discretized in terms of the stream function variable $\psi$ on the spatial domain $[0,1]^2$ using a $p$-order central finite difference discretization, which by default is $p=2$. 
The constants are taken to be $F=1600$, $\varepsilon=10^{-5}$, and $ A=2\times10^{-11}$, with $A$ being the most variable, and potentially an order of magnitude smaller in some instances.
Homogeneous Dirichlet boundary conditions are applied, and therefore the zero-valued boundary points are not explicitly stored in the computational state-space variable. The dimension is taken to be $n=127^2=16,129$ by default.
The value of $J$ is calculated by the Arakawa approximation~\citep{arakawa1966computational,ferguson2008numerical}. The operator $\Delta^3$ is computed offline as the cube of the sparse Laplacian.

We provide several different ways of solving the Helmhotlz equation in \cref{eq:QGSO}, which by default is solved through an offline pivoted sparse Cholesky decomposition. A multigrid solver and the ability to solve via \texttt{GMRES} are also provided.

Additionally, as several practical Data Assimilation applications require the computational of physical distance between state-space variables. We define the distance function as
\begin{align*}
  d(i, j) = \sqrt{{(i_x - j_x)}^2 + {(i_y - j_y)}^2},
\end{align*}
where $(i_x, i_y)$ and $(j_x, j_y)$ are the $i$ and $j$ components of two distinct state space variables when the state space is transformed into a two dimensional grid.

This particular implementation of QGSO was written specifically for~\citep{popov2018bayesian}.

The problem also has both Jacobian vector products and Adjoint vector products configured to work with the same linear system solvers as that of the right hand side. A Jacobian approximation is also provided. 

\begin{figure}
\begin{center}
  \includegraphics[width=0.90\textwidth]{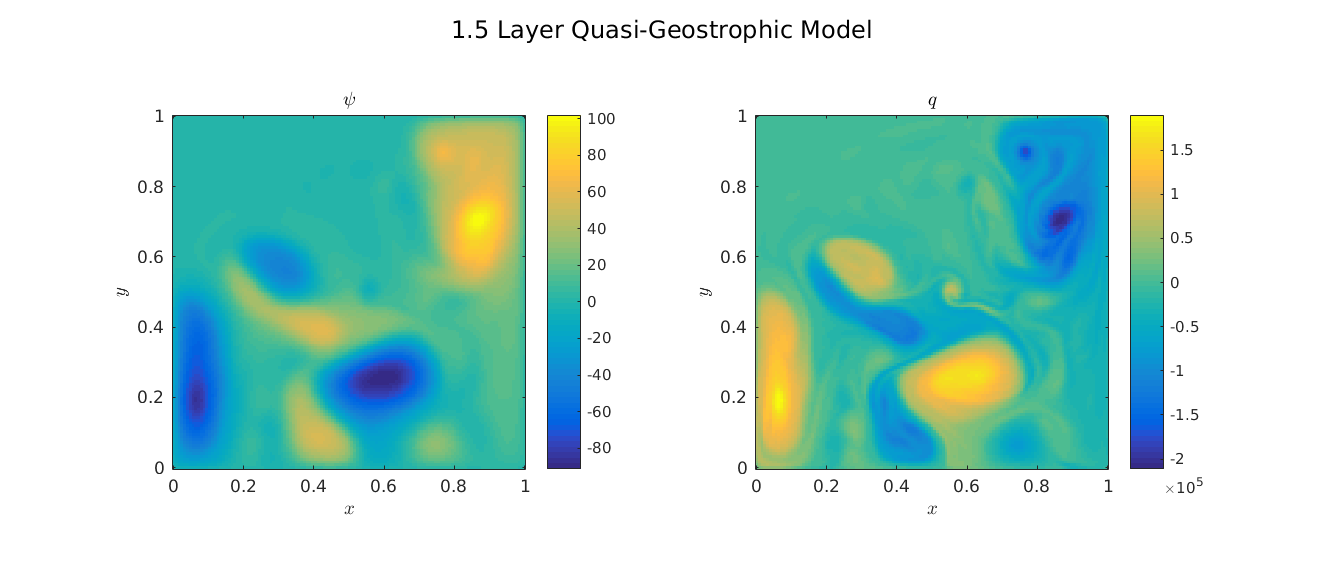}%
\end{center}
\caption{A typical model state of the 1.5 layer quasi-geostrophic model (QGSO). The left graph represents the stream function while the right represents vorticity.}
\end{figure}

\getpresettable{qgso}
\problemdescription{Gray--Scott}{grayscott}{$n^2$ ($128^2$)}

The Gray--Scott model \cite{gray1983autocatalytic} is a reaction--diffusion problem that simulates the chemical reaction
\begin{align*}
    U + 2V &\to 3V, \\
    V &\to P.
\end{align*}
The corresponding PDEs are given by
\begin{align}
    \begin{split}
        \frac{\partial u}{\partial t} &= \varepsilon_1 \Delta u - u v^2 + f(1 - u), \\
        \frac{\partial v}{\partial t} &= \varepsilon_2 \Delta v + u v^2 - (f + k) v, \\
    \end{split}
\end{align}
with periodic boundary conditions.  $\varepsilon_1$ and $\varepsilon_2$ are diffusion rates and $f$ and $k$ are reaction rates.  This problem is discretized on a uniform 2D grid with second order finite differences.

\getpresettable{grayscott}

\section{Conclusions}

We provide a new \MATLAB{}-based framework of the creation, manipulation, and usage of problems posed as first order ordinary differential equations.  OTP has an extensible, easy-to-use interface.  In the future, additional problems and functionality will be added.  Of special interest are more PDEs and large-scale problems.

\begin{acknowledgments}
 Special thanks to Mahesh Narayanamurthi, S. Ross Glandon, Arash Sarshar, Bibek Regmi, and the rest of the Computational Science Lab at Virginia Tech for their patience and support of this project.
\end{acknowledgments}

\bibliographystyle{plain}
\bibliography{main}

\end{document}